\documentclass{article}

\usepackage{amsthm}
\usepackage{amsfonts}
\usepackage{amssymb}
\usepackage{amsgen}
\usepackage{amsmath}
\usepackage{amsopn}
\usepackage{verbatim}
\usepackage{xypic}
\usepackage{multicol}
\usepackage{eepic}
\usepackage{upref}
\usepackage{pgf}
\usepackage{tikz}
\usetikzlibrary{arrows}

\theoremstyle{plain}
\newtheorem{thm}{Theorem}

\newtheorem{prop}{Proposition}

\newtheorem{conj}{Conjecture}

\theoremstyle{definition}

\theoremstyle{remark}
\newtheorem{rem}{Remark}

 \font\cyr=wncyr10
 \newcommand{\nc}{\newcommand}

\DeclareMathOperator{\Hom}{{Hom}}
\DeclareMathOperator{\Gr}{{Gr}}

\DeclareMathOperator{\dch}{{dch}}
\DeclareMathOperator{\Lie}{{Lie}}
\DeclareMathOperator{\MT}{{MT}}

 \nc{\wht}{{\widehat}}
 \nc{\bwg}{{\bigwedge}}
 \nc{\mmu}{{\boldsymbol{\mu}}}
 \nc{\mal}{{{\scriptstyle \maltese}}}
 \nc{\fA}{{\mathfrak A}}
 \nc{\HH}{{\mathfrak H}}
 \nc{\ra}{\rightarrow}
 \nc{\ors}{{\vec s\,}}
 \nc{\os}{{\overset}}
 \nc{\G}{{\mathbb G}}
 \nc{\Z}{{\mathbb Z}}
 \nc{\R}{{\mathbb R}}
 \nc{\N}{{\mathbb N}}
 \nc{\ZN}{{\mathbb Z_{\ge 0}}}
 \nc{\Q}{{\mathbb Q}}
 \nc{\C}{{\mathbb C}}
 \nc{\MPV}{{\mathcal{MPV}}}
 \renewcommand{\P}{{\mathbb P}}
 \nc{\tB}{{\tilde B}}
 \nc{\Li}{{\rm Li}}
 \nc{\suf}{{\ast\,}}
 \nc{\sufq}{{\ast_q\,}}
 \nc{\gam}{{\gamma}}
 \nc{\gG}{{\Gamma}}
 \nc{\om}{{\omega}}
 \nc{\vep}{{\varepsilon}}
 \nc{\ga}{{\alpha}}
 \nc{\gl}{{\lambda}}
 \nc{\gb}{{\beta}}
 \nc{\gd}{{\delta}}
 \nc{\gs}{{\sigma}}
 \nc{\gS}{{\Sigma}}
 \nc{\gO}{{\Omega}}
 \nc{\sif}{{\mathcal S}}
 \nc{\gt}{{\tau}}
 \nc{\Lra}{\Longrightarrow}
 \nc{\lra}{\longrightarrow}
 \nc{\fS}{{\mathfrak S}}
 \nc{\DD}{{\mathfrak D}}

 \nc{\Llra}{\Longleftrightarrow}
 \nc{\ol}{\overline}
 \nc{\lms}{\longmapsto}
 \nc{\cv}{{{\mathsf c}{\mathsf v}}}
 \nc{\zq}{{\zeta_q}}
 \nc\qup{{q\uparrow 1}}
 \nc{\us}{\underset}
 \nc{\tn}{{\tilde{n}}}
 \nc{\gD}{{\Delta}}
 \nc{\bk}{{\bf k}}
 \nc{\bi}{{\bf i}}
 \nc{\bfone}{{\bf 1}}
 \nc{\bfs}{{\bf s}}
 \nc{\bfx}{{\bf x}}
 \nc{\bfY}{{\bf Y}}
 \nc{\QX}{{\Q\langle \bfX\rangle}}
 \nc{\QY}{{\Q\langle \bfY\rangle}}
 \nc{\CX}{{\C\langle \bfX\rangle}}
 \nc{\CY}{{\C\langle \bfY\rangle}}
 \nc{\QXX}{{\Q\langle\!\langle \bfX\rangle\!\rangle}}
 \nc{\QYY}{{\Q\langle\!\langle \bfY\rangle\!\rangle}}
 \nc{\CXX}{{\C\langle\!\langle \bfX\rangle\!\rangle}}
 \nc{\CYY}{{\C\langle\!\langle \bfY\rangle\!\rangle}}

 \nc{\bbA}{{\mathbb A}}
 \nc{\bbB}{{\mathbb B}}
 \nc{\bbC}{{\mathbb C}}
 \nc{\bbD}{{\mathbb D}}
 \nc{\bbE}{{\mathbb E}}
 \nc{\bbF}{{\mathbb F}}
 \nc{\bbG}{{\mathbb G}}
 \nc{\bbH}{{\mathbb H}}
 \nc{\bbI}{{\mathbb I}}
 \nc{\bbJ}{{\mathbb J}}
 \nc{\bbK}{{\mathbb K}}
 \nc{\bbL}{{\mathbb L}}
 \nc{\bbM}{{\mathbb M}}
 \nc{\bbN}{{\mathbb N}}
 \nc{\bbO}{{\mathbb O}}
 \nc{\bbP}{{\mathbb P}}
 \nc{\bbQ}{{\mathbb Q}}
 \nc{\bbR}{{\mathbb R}}
 \nc{\bbS}{{\mathbb S}}
 \nc{\bbT}{{\mathbb T}}
 \nc{\bbU}{{\mathbb U}}
 \nc{\bbV}{{\mathbb V}}
 \nc{\bbW}{{\mathbb W}}
 \nc{\bbX}{{\mathbb X}}
 \nc{\bbY}{{\mathbb Y}}
 \nc{\bbZ}{{\mathbb Z}}
 \nc{\bba}{{\mathbb a}}
 \nc{\bbb}{{\mathbb b}}
 \nc{\bbc}{{\mathbb c}}
 \nc{\bbd}{{\mathbb d}}
 \nc{\bbe}{{\mathbb e}}
 \nc{\bbf}{{\mathbb f}}
 \nc{\bbg}{{\mathbb g}}
 \nc{\bbh}{{\mathbb h}}
 \nc{\bbi}{{\mathbb i}}
 \nc{\bbk}{{\mathbb k}}
 \nc{\bbl}{{\mathbb l}}
 \nc{\bbm}{{\mathbb m}}
 \nc{\bbn}{{\mathbb n}}
 \nc{\bbo}{{\mathbb o}}
 \nc{\bbp}{{\mathbb p}}
 \nc{\bbq}{{\mathbb q}}
 \nc{\bbr}{{\mathbb r}}
 \nc{\bbs}{{\mathbb s}}
 \nc{\bbt}{{\mathbb t}}
 \nc{\bbu}{{\mathbb u}}
 \nc{\bbv}{{\mathbb v}}
 \nc{\bbw}{{\mathbb w}}
 \nc{\bbx}{{\mathbb x}}
 \nc{\bby}{{\mathbb y}}
 \nc{\bbz}{{\mathbb z}}

 \nc{\calA}{{\mathcal A}}
 \nc{\calB}{{\mathcal B}}
 \nc{\calC}{{\mathcal C}}
 \nc{\calD}{{\mathcal D}}
 \nc{\calE}{{\mathcal E}}
 \nc{\calF}{{\mathcal F}}
 \nc{\calG}{{\mathcal G}}
 \nc{\calH}{{\mathcal H}}
 \nc{\calI}{{\mathcal I}}
 \nc{\calJ}{{\mathcal J}}
 \nc{\calK}{{\mathcal K}}
 \nc{\calL}{{\mathcal L}}
 \nc{\calM}{{\mathcal M}}
 \nc{\calN}{{\mathcal N}}
 \nc{\calO}{{\mathcal O}}
 \nc{\calP}{{\mathcal P}}
 \nc{\calQ}{{\mathcal Q}}
 \nc{\calR}{{\mathcal R}}
 \nc{\calS}{{\mathcal S}}
 \nc{\calT}{{\mathcal T}}
 \nc{\calU}{{\mathcal U}}
 \nc{\calV}{{\mathcal V}}
 \nc{\calW}{{\mathcal W}}
 \nc{\calX}{{\mathcal X}}
 \nc{\calY}{{\mathcal Y}}
 \nc{\calZ}{{\mathcal Z}}
  \nc{\cala}{{\mathcal a}}
 \nc{\calb}{{\mathcal b}}
 \nc{\calc}{{\mathcal c}}
 \nc{\cald}{{\mathcal d}}
 \nc{\cale}{{\mathcal e}}
 \nc{\calf}{{\mathcal f}}
 \nc{\calg}{{\mathcal g}}
 \nc{\calh}{{\mathcal h}}
 \nc{\cali}{{\mathcal i}}
 \nc{\calj}{{\mathcal j}}
 \nc{\calk}{{\mathcal k}}
 \nc{\call}{{\mathcal l}}
 \nc{\calm}{{\mathcal m}}
 \nc{\caln}{{\mathcal n}}
 \nc{\calo}{{\mathcal o}}
 \nc{\calp}{{\mathcal p}}
 \nc{\calq}{{\mathcal q}}
 \nc{\calr}{{\mathcal r}}
 \nc{\cals}{{\mathcal s}}
 \nc{\calt}{{\mathcal t}}
 \nc{\calu}{{\mathcal u}}
 \nc{\calv}{{\mathcal v}}
 \nc{\calw}{{\mathcal w}}
 \nc{\calx}{{\mathcal x}}
 \nc{\caly}{{\mathcal y}}
 \nc{\calz}{{\mathcal z}}

 \nc{\frakA}{{\mathfrak A}}
 \nc{\frakB}{{\mathfrak B}}
 \nc{\frakC}{{\mathfrak C}}
 \nc{\frakD}{{\mathfrak D}}
 \nc{\frakE}{{\mathfrak E}}
 \nc{\frakF}{{\mathfrak F}}
 \nc{\frakG}{{\mathfrak G}}
 \nc{\frakH}{{\mathfrak H}}
 \nc{\frakI}{{\mathfrak I}}
 \nc{\frakJ}{{\mathfrak J}}
 \nc{\frakK}{{\mathfrak K}}
 \nc{\frakL}{{\mathfrak L}}
 \nc{\frakM}{{\mathfrak M}}
 \nc{\frakN}{{\mathfrak N}}
 \nc{\frakO}{{\mathfrak O}}
 \nc{\frakP}{{\mathfrak P}}
 \nc{\frakQ}{{\mathfrak Q}}
 \nc{\frakR}{{\mathfrak R}}
 \nc{\frakS}{{\mathfrak S}}
 \nc{\frakT}{{\mathfrak T}}
 \nc{\frakU}{{\mathfrak U}}
 \nc{\frakV}{{\mathfrak V}}
 \nc{\frakW}{{\mathfrak W}}
 \nc{\frakX}{{\mathfrak X}}
 \nc{\frakY}{{\mathfrak Y}}
 \nc{\frakZ}{{\mathfrak Z}}
 \nc{\fraka}{{\mathfrak a}}
 \nc{\frakb}{{\mathfrak b}}
 \nc{\frakc}{{\mathfrak c}}
 \nc{\frakd}{{\mathfrak d}}
 \nc{\frake}{{\mathfrak e}}
 \nc{\frakf}{{\mathfrak f}}
 \nc{\frakg}{{\mathfrak g}}
 \nc{\frakh}{{\mathfrak h}}
 \nc{\fraki}{{\mathfrak i}}
 \nc{\frakj}{{\mathfrak j}}
 \nc{\frakk}{{\mathfrak k}}
 \nc{\frakl}{{\mathfrak l}}
 \nc{\frakm}{{\mathfrak m}}
 \nc{\frakn}{{\mathfrak n}}
 \nc{\frako}{{\mathfrak o}}
 \nc{\frakp}{{\mathfrak p}}
 \nc{\frakq}{{\mathfrak q}}
 \nc{\frakr}{{\mathfrak r}}
 \nc{\fraks}{{\mathfrak s}}
 \nc{\frakt}{{\mathfrak t}}
 \nc{\fraku}{{\mathfrak u}}
 \nc{\frakv}{{\mathfrak v}}
 \nc{\frakw}{{\mathfrak w}}
 \nc{\frakx}{{\mathfrak x}}
 \nc{\fraky}{{\mathfrak y}}
 \nc{\frakz}{{\mathfrak z}}

 \nc{\sha}{{\mbox{\cyr x}}}

\begin{document}

\title{Multiple polylogarithm values at roots of unity}
\author{Jianqiang Zhao}
\date{}
\maketitle
\begin{center}
{Department of Mathematics, Eckerd College, St. Petersburg, FL 33711, USA}
\end{center}

\medskip

\noindent{\bf Abstract.} For any positive integer $N$ let $\mmu_N$ be the
group of the $N$th roots of unity. In this note we shall study the $\Q$-linear relations
among values of multiple polylogarithms evaluated at $\mmu_N$. We show that
the standard relations considered by Racinet do not provide all the possible
relations in the following cases: (i) level $N=4$, weight $w=3$ or $4$, and (ii) $w=2$,
$7<N<50$, and $N$ is a power of 2 or 3, or $N$ has at least two prime factors.
We further find some (presumably all) of the missing relations in (i) by using the
octahedral symmetry of $\P^1-(\{0,\infty\}\cup \mmu_4)$. We also prove some
other results when $N=p$ or $N=p^2$ ($p$ prime $\ge 5$) by using the
motivic fundamental group of $\P^1-(\{0,\infty\}\cup\mmu_N)$.

\medskip
\section{Introduction}

Double shuffle relations have played significant roles in the study of
multiple zeta values in recent years. These relations can be easily
generalized to the special values of
multiple polylogarithms at roots of unity (MPV for short):
\begin{equation}\label{equ:seriesdefn}
Li_{s_1,\dots,s_n}(\zeta_1,\dots,\zeta_n):=\sum_{k_1>\dots>k_n>0}
\frac{\zeta_1^{k_1} \cdots \zeta_n^{k_n}}{k_1^{s_1}\cdots k_n^{s_n}},\quad(s_1,\zeta_1)\ne (1,1),
\end{equation}
where $\zeta_j$ are all $N$th roots of unity for $j=1,\dots,n$. We call $N$ the
\emph{level} and $w:=s_1+\cdots+s_n$ the \emph{weight}. One of our major
interests is to find the dimension $d(w,N)$
of the $\Q$-vector space $\MPV(w,N)$ spanned by these values.

In \cite{Rac}, Racinet listed the following relations:
double shuffle relations and their regularized versions (by regularizing
both the divergent integrals and divergent series representing MPVs),
(regularized) distributions, and weight one relations.
We will call these \emph{standard relations}.
By intensive MAPLE computation it is shown \cite{Zpolyrel} that in
many cases these relations cannot produce all the possible $\Q$-linear
relations. For example, standard relations imply only $d(3,4)\le 9$ which is
one more than the bound given by \cite[5.25]{DG}.
Concretely, from standard relations MAPLE confirms the following
\begin{equation}\label{fact}
\aligned
 \text{Fact: }&\text{Every MPV of weight 3 and level 4 can be written explicitly} \\
 &\text{as a linear combination of the nine MPVs appearing in \eqref{equ:conj},}
\endaligned
\tag{$\ast$}
\end{equation}
and no further $\Q$-linear relations between these values can be deduced from
the standard relations. But by GiNac \cite{GiNac} and EZface \cite{EZface}
the following is found numerically (see \cite[Remark 10.1]{Zpolyrel})
\begin{multline}\label{equ:conj}
5Li_{1,2}(-1,-i)=46Li_{1,1,1}(i,1,1)-7Li_{1,1,1}(-1,-1,i)-13Li_{1,1,1}(i,i,i) +13Li_{1,2}(-i,i)\\
-Li_{1,1,1}(-i,-1,1)+25Li_{1,1,1}(-i,1,1)-8Li_{1,1,1}(i,i,-1) +18Li_{2,1}(-i,1).
\end{multline}
In this note, we shall prove
\eqref{equ:conj} from the explicit relations in \eqref{fact}
by using the octahedral symmetry of
$\P^1-(\{0,\infty\}\cup \mmu_4)$ where $\mmu_N$ denotes the set of $N$th roots of unity.
We can treat the weight 4 and level 4 case in a similar fashion.

The result above in level four case shows that the standard relations are not
always sufficient to determine $d(w,N)$.
On the other hand, Deligne and Goncharov \cite{DG}
find some closed formulae for upper
bounds of $d(w,N)$ for all $w$ and $N$ by studying the motivic fundamental group of
$\P^1-(\{0,\infty\}\cup \mmu_N)$ in the framework of mixed Tate motives over
the ring of $S$-integers of a number field.
\begin{prop} \label{prop:DG}
\emph{(\cite[5.25]{DG})} Let $N$ be a positive integer.
Let $\varphi$ be the Euler's totient function
and $\nu(N)$ be the number of prime factors of $N$.
Then $d(w,N)\le D(w,N)$
where $D(w,N)$ are defined by the formal power series
$$1+\sum_{w=1}^\infty D(w,N)t^w=
\left\{
    \begin{array}{ll}
      (1-t^2-t^3)^{-1}, & \hbox{if $N=1$;} \\
      (1-t-t^2)^{-1}, & \hbox{if $N=2$;} \\
     \big(1-\big(\frac{\varphi(N)}2+\nu(N)\big)t+ \big(\nu(N)-1 \big)t^2 \big)^{-1}, & \hbox{if $N\ge 3$.} \\
    \end{array}
  \right.
$$
\end{prop}
As an example one can find the bound $d(w,4)\le 2^w$ for all weight $w$ and
by a variant of a conjecture of Grothendieck it can
be deduced that for $N=4$ equality should always hold.
Further, Goncharov \cite{DG} has defined some pro-nilpotent
Lie algebra $D$ graded by weight and depth using standard
relations and then related $D$ to the motivic fundamental Lie algebra
in lower depth cases. Using this we shall improve the bounds in the
above proposition for all levels $N=p$ and $N=p^2$ ($p\ge 5$ a prime)
in the last section. From this improved bound we find that when
$7<N<50$ standard relations are incomplete if $N$ is a power
of 2 or 3, or $N$ has at least two prime factors, though we
don't know how to find missing relations explicitly.

\section{The motivic setup}
Let us review some constructions from \cite{DG}.
Fix a positive integer $N$. Let $k$ be the field over $\Q$ generated
by a primitive $N$th root of unity and $\calO$ its ring of integers. Let
$\Gamma$ be the $\Q$-vector subspace of $\calO[1/N]^*\otimes\Q$ generated by
$1-\zeta$ for $\zeta$ an $N$th root of unity other than 1.
We refer to \cite{DG} for the definition of the tannakian
category of mixed Tate motives $\MT(k)_\Gamma$,
the definition of the fibre functor
$\om(M):=\bigoplus_n\Hom(\Q(n),\Gr^W_{-2n}(M))$,
and of the corresponding motivic fundamental group $G$.
Let $\calL$ be the Lie algebra generated
by generators $e_x,x\in\{0,\infty\}\cup \mmu_N,$ with the only
relation $\sum e_x=0$. Set the degree of $e_x$ to be 1. Then
$\calL$ is freely generated by the $e_x,x\ne \infty$. Let $\Pi$
be the pro-nilpotent group
$ \lim_n \exp(\calL/\text{deg}^{>n})$. The completed enveloping algebra
of $\calL$ is
$\Q\langle\!\langle e_\zeta \rangle\!\rangle_{\zeta\in\mmu_N\cup \{0\}}$
with the coproduct $\gD: e_\zeta\mapsto e_\zeta\otimes 1+1\otimes e_\zeta$,
and $\Pi$ is the subscheme of the group-like elements under $\gD$.

Let $H_\om$ be the group of $\mmu_N$-equivariant automorphisms of $\Pi$
defined in \cite[\S5.8]{DG}.
Set $G_\om:=\om(G),$ $U_\om:=\ker(G_\om\to \G_m)$ and $V_\om:=\ker(H_\om\to \G_m)$.
There are morphisms
of group schemes $G_\om\to H_\om$ and $U_\om\to V_\om$
compatible with the semi-direct product structure
\begin{equation}\label{equ:morphs}
    \iota: G_\om =\G_m\ltimes U_\om\to H_\om =\G_m\ltimes V_\om.
\end{equation}

Fix an embedding $\gs:k\to \C$. Let
$\Omega=\sum_{a\in \mmu_N\cup\{0\}} \frac{dz}{z-\gs(a)}e_a.$
Paths $\frakp$ in $\C^*-\mmu_N$ provide complex points
$I(\frakp)=\sum_{n=0}^\infty  \int_{\frakp} \Omega^{\circ n}$ in $\Pi$
with coordinates given by Chen's iterated integrals so that
$I(\frakp\text{ followed by } \frakq)=I(\frakq)I(\frakp).$
At punctured points $x\in \{0,\infty\}\cup\mmu_N$ one may use
the tangential base points to extend the above construction
to paths beginning or ending at $x$. This is a
regularization procedure of divergent integrals with the
counter terms prescribed by the tangent vectors. For example, one may
define the straight line path from $1_0$ (tangent vector 1 at 0)
to $(-1)_1$ (see \cite[p.~85]{D1}). This
gives rise to $\dch(\gs)$ whose coefficient of
$e_0^{s_1-1}e_{\zeta_1}\dots e_0^{s_n-1}e_{\zeta_n}$
is $(-1)^nLi_{s_1,\dots,s_n}(\gs(1/\zeta_1),
\gs(\zeta_1/\zeta_2),\dots,\gs(\zeta_{n-1}/\zeta_n))$
if $(s_1,\zeta_1)\ne (1,1)$. Note that
$\dch(\gs)$ is the unique group-like element in
$\Q\langle\!\langle e_\zeta \rangle\!\rangle_{\zeta\in\mmu_N\cup \{0\}}$
whose coefficients
of $e_0$ and $e_1$ are both 0 with that property.

\section{Standard relations in $\MPV(w,4)$}\label{sec:srel}
Fix $N=4$ throughout this section.
In what follows we shall describe a process suggested by Deligne
to verify that the standard relations
cannot produce all possible $\Q$-linear relations in $\MPV(w,4)$.

Observe that the vector space freely generated by basis vectors corresponding
to regularized MPVs of weight 3 and level 4 is dual to the degree 3
part of the free associative algebra
$\Q\langle\!\langle e_0,e_\zeta \rangle\!\rangle_{\zeta\in\mmu_4}$.
One expects that $\MPV(w,4)_{w\ge 1}$ is a weighted polynomial algebra with
the following number
of generators in low weights: 2,1, 2, and 3 (in weight 1 to 4).
The reason is that $\iota(\Lie U_\om)$ should be a Lie algebra freely
generated by one element in each degree so that the dimension of degree $n$
part is given by $\frac1n\sum_{d|n}\mu(n/d) 2^d-\delta_{1,n}$
(see \cite[Ch.\ II, \S3 Thm.\ 2]{Bourbaki}) which is the
number of binary Lyndon words of length $n$ when $n>1$.

If one takes the space $\MPV(3,4)$ modulo the products of MPVs in lower weights one
should get a quotient space of dimension 2 and one knows its dimension $\le 2$. However,
one can only prove the bound 3 by the standard relations as follows.
Consider the subspace generated by the following elements:

\medskip
I. For each MPV of weight 1 and MPV of weight 2, the linear combination
of weight 3 MPVs expressing their shuffle product by Chen's iterated integrals.
This gives $5\times 25$ elements.

\medskip
II. Same for the stuffle product (quasi-shuffle as called by Hoffman)
corresponding to the coproduct $\gD_*$ in \cite{Rac}.
Plainly, these elements are linear combination of weight 3 MPVs expressing
their stuffle products by series expansions \eqref{equ:seriesdefn}.
This gives $4\times 20$ elements.

\medskip
III. Distribution relations in weight 3: expressing the coefficient
of a ``convergent'' monomial in $e_0$ and $e_\ga$ ($\ga=\pm1$) as a multiple
of the sum of the coefficients of the monomials deduced from it by replacing
$e_\ga$ by $e_\gb$ with $\gb^2=\ga$. There are 12 ``convergent'' monomials
to consider:
$(e_0\text{ or }e_{-1})(e_0\text{ or }e_1\text{ or }e_{-1})(e_1\text{ or }e_{-1}).$

\medskip
IV. Regularized distribution relations in weight 3 (see \cite[Prop.~2.26]{Rac},
change $\gs^{n/d}=1$ to $\gs^d=1$). There are 6 of these.

\medskip
All of these elements from I to IV can be put together to form a $223\times 125$ matrix.
By MAPLE one can verify that it only has rank 122. Moreover, it is not hard to find three
linearly independent relations over $\Q$ which produce three vectors
in the degree 2 part of the Lie algebra denoted by $\frak{dmrd}_0$ in \cite{Rac}.
In weight 4 the same procedure produces eight vectors instead of three which means one needs
five more relations besides the standard ones.

\section{Octahedral symmetry of $\P^1-(\{0,\infty\}\cup \mmu_4)$}
The Riemann sphere $\P^1$ punctured at $\{0,\infty,\pm1,\pm i\}$ clearly possesses
an octahedral symmetry. However, one should use a system of tangent vectors
stable (up to multiplication by roots of unity) by the octahedral group to
use this symmetry: one should use the straight line path from $1_0$ to $(-2)_1$.
Define the $\C$-linear map
$\rho: \Q\langle\!\langle e_0,e_\zeta \rangle\!\rangle_{\zeta\in\mmu_4}
\lra\Q\langle\!\langle e_0,e_\zeta \rangle\!\rangle_{\zeta\in\mmu_4}$ such that
$$\rho: e_0\to e_1\to e_i\to e_0, \qquad e_\infty\to e_{-1}\to e_{-i}\to e_\infty.$$
\begin{multicols}{2}
Let $0<\vep<1/3$ and $C_\vep$ be the path $A_1\dots A_6$ in the complex plane
shown in the right picture, where $A_2, A_4$ and $A_6$ are the quarter circles
of radii $\vep$, $2\vep$ and $2\vep$, respectively, oriented clockwise.
Then by the property of iterated integrals
\begin{center}
\begin{tikzpicture}[scale=0.7pt]
\draw[gray,very thin](-0.5,0) -- (3.1,0);
\draw[gray,very thin] (0,3.1) -- (0,-0.5);
\draw[->](1.4,0) -- (1.5,0) node [anchor=north,below=0.5pt] {$A_1$};
\draw[->](0,1.6) -- (0,1.5) node [left=0.2pt] {$A_5$};
\draw (-0.2,0) node [below=1pt] {$0$};
\draw (0,0) node [rotate=-45,above=3pt] {$A_6$};
\draw (0,0) node {$\bullet$};
\draw (3,0) node {$\bullet$};
\draw (0,3) node {$\bullet$};
\draw (0,0) node[rotate=45,right=45pt]  {$A_3$};
\draw (0.2,0)  -- (2.6,0);
\draw (3,0) node[rotate=45,above=5pt] {$A_2$};
\draw (0,3) node[rotate=45,below=5pt]  {$A_4$};
\draw (-0.2,3) node {$i$};
\draw [->](2.6,0) arc (180:135:0.4) ;
\draw (3,0.4) arc (90:180:0.4);
\draw (3,0) node[below=1pt] {$1$};
\draw [->](3,0.4) arc (0:45:2.6);
\draw (3,0.4) arc (0:90:2.6);
\draw [->](0.4,3) arc (0:-45:0.4);
\draw (0,2.6) arc (270:360:0.4);
\draw (0,2.6) -- (0,0.2);
\draw [->](0,0.2) arc (90:45:0.2);
\draw (0.2,0) arc (0:90:0.2);
\end{tikzpicture}
\end{center}
\end{multicols}

\vskip-.7cm

\begin{equation}\label{equ:loopintegral}
 1=\int_{C_\vep} \sum_{n=0}^\infty \Omega^{\circ n}=\int_{A_6}  \sum_{n=0}^\infty \Omega^{\circ n}
\int_{A_5}  \sum_{n=0}^\infty \Omega^{\circ n}\cdots \int_{A_1}  \sum_{n=0}^\infty \Omega^{\circ n}.
 \end{equation}
Replacing the straight path $1_0\to (-1)_1$ by
the straight path $1_0\to (-2)_1$ changes $\dch(\gs)$ to
$\dch'(\gs)=\exp\big((\log 2) e_1\big)\dch(\gs)$.
So regularized integral over the path $A_1$ followed by $A_2$ gives
$I=\exp(-2\pi i e_1/4)\dch'(\gs)$.
Thus \eqref{equ:loopintegral} yields
$\rho^2(I) \rho(I) I=1$, i.e.,
\begin{equation}  \label{equ:octa}
\exp(-2Li_1(i) e_0)\rho^2(\dch(\gs))\exp(-2Li_1(i)e_i)
\rho(\dch(\gs)) \exp\left(-2Li_1(i)e_1\right) \dch(\gs) =1.
\end{equation}
By comparing the coefficient of $e_2e_1^2$ in \eqref{equ:octa} and using Fact \eqref{fact}
one finally arrives at \eqref{equ:conj}.

We now consider the Lie algebra $\calG:=\iota(\Lie U_\om)$ which is the image of
the map in \eqref{equ:morphs}, equipped with Ihara's Lie bracket $\{\ ,\ \}$
(see \cite[(5.13.6)]{DG}). One can use the additional
octahedral relation to find the correct number of generators
of $\calG$ in degrees up to 3:
\begin{align*}
    \text{deg 1: }&\ v_1=2e_{-1}+2e_1+e_{-i}+e_i\\
    \text{deg 2: }&\  v_2=[e_0,e_i]-[e_0,e_{-i}]+[e_1,e_i]+[e_{-i},e_1]+[e_{-i},e_i]\\
    \text{deg 3: }&\  \{v_1,v_2\} \text{ and }  v_3,
\end{align*}
where $v_3$ is obtained by MAPLE and is too long to write down here.
Define on the Lie algebra an involution $\gs: e_0 \leftrightarrow e_1, e_i
\leftrightarrow e_{-i}, e_{-1} \leftrightarrow e_\infty$ in the
octahedral group.
Then one can check by MAPLE that $v+\gs(v)=0$ for all the above generators.

Next we look at weight 4 case. We find that the upper bound produced
by the standard relations is only $d(4,4)\le 21$ \cite[Table 2]{Zpolyrel}. By exactly the
same approach as in weight 3 one can verify that all the five missing
relations can be produced by applying the octahedral symmetry. For the
five octahedral relations one can use the vanishing of the coefficients of
$e_{-i}e_0^2e_{-i},\  e_{-i}e_0^2e_{-1},\  e_{-i}e_0^2e_{i},\  e_{-i}e_0^2e_{1}$,
and $(e_{-i}e_0)^2$ in \eqref{equ:octa}.
Moreover, a process similar to the one sketched in \S\ref{sec:srel}
yields exactly three generators in the space
of the degree 4 part of the Lie algebra $\calG$:
$$\text{deg 4: }\ \{v_1,v_3\},\quad \{v_1,\{v_1,v_2\}\}, \text{ and }v_4$$
where $v_4$ is a complicated new vector obtained by MAPLE. This verifies
up to degree 4 the following fact: $\calG$ is a Lie algebra freely
generated by one generator in each degree.

\begin{rem} For arbitrary $N$ the set $\{0,\infty\}\cup \mmu_N$
possesses the dihedral symmetry. However, when we apply this
in $N=4$ cases (weight 3 or 4) all the relations produced are
already known to us by standard relations. This should be true
in general (cf. \cite[Thm.~4.1]{G1}).
\end{rem}

\section{Level $p$ and $p^2$ cases ($p\ge 5$ a prime)}\label{sec:pp}
Let $p\ge 5$ be prime throughout this section. From Prop.~\ref{prop:DG}
we see that $d(2,p)\le (p+1)^2/4$. We will improve this bound
in Theorem~\ref{thm:p}.

Recall that $\calG=\iota(\Lie U_\om)$.
By \cite[6.13]{DG} one may safely replace $\calG_N^{(\ell)}$ by
$\calG$ throughout \cite{G1}. In particular, Goncharov's results
can now be used to study the structure of Galois Lie algebra
$\calG_{\bullet,\bullet}(\mmu_N)=\sum_{w\ge 1,l\ge 1}\calG_{-w,-l}(\mmu_N)$
graded by the weight $w$ and depth $l$. The weight filtration
is defined via the action of $\G_m$ and the depth filtration  comes from
the descending central series of $\ker f_*$, where
$f_*:\pi_1(\G_m-\mmu_N,x) \lra \pi_1(\G_m)$
for any base point $x$ is induced from the embedding
$\G_m-\mmu_N\hookrightarrow\G_m$.

\begin{thm}\label{thm:p} Let $p\ge 5$ be a prime. Then
\begin{equation}\label{equ:pbdd}
 d(2,p)\le \frac{(5p+7)(p+1)}{24}.
\end{equation}
If Grothendieck's period conjecture (see \cite[5.27(c)]{DG})  is true then
the equality holds and all the $\Q$-linear relations in $\MPV(2,p)$
follow from the standard relations.
\end{thm}
\begin{proof}
Let $\{\ ,\ \}$ be the Lie bracket in $\calG_{\bullet,\bullet}(\mmu_N)$.
Now in depth 1 one has $\dim \calG_{-2,-1}(\mmu_p)=(p-1)/2$
(see \cite[Theorem~2.1]{G1}). By \cite[Cor.~2.16]{G1}
\begin{equation}\label{equ:diag2}
 \dim \calG_{-2,-2}(\mmu_p)=\frac{(p-1)(p-5)}{12}.
\end{equation}
For any positive integer $N$ define
\begin{equation}
 \label{equ:beta}
 \beta_N:\bigwedge^2\calG_{-2,-1}(\mmu_N)  \lra \calG_{-2}(\mmu_N), \quad
  a\wedge b  \lms\{a,b\}.
\end{equation}
Then
\begin{equation}\label{equ:kerbeta}
 \dim(\ker\beta_p)=\frac12\frac{p-1}{2}\left(\frac{p-1}2-1\right)-\frac{(p-1)(p-5)}{12}
=\frac{p^2-1}{24}.
\end{equation}

Let $SR$ be the affine subgroup of $\Pi$ defined by only those
polynomial equations satisfied by the coefficients of $\dch(\gs)$
which are deduced from the standard relations, plus ``$2\pi i=0$'',
as explained in \cite[5.22]{DG}. This is the group $\mathsf{DMRD}_0$
studied by Racinet in \cite[\S3.2, Thm.~I]{Rac}. Its Lie algebra $\Lie SR$ is
graded by weight and depth, independent of the embedding
by Prop.~4.1 of op. cit. and contains $\calG$ by \cite[5.22]{DG}.
Further, Goncharov shows \cite[\S7.7]{G1} that the
standard relations provide a complete list
of constraints on the diagonal part of the Lie algebra $\calG$ in depth $\le 2$,
yielding $(\Lie SR)_{m,m}=\calG_{m,m}$ for $m=-1,-2$.
Together with \eqref{equ:kerbeta} this means in the proof of
\cite[5.25]{DG} one can decrease
the bound $D(2,p)$ by $(p^2-1)/24$ and arrive at the bound
$(5p+7)(p+1)/24$.

Grothendieck's conjecture implies that
$\calG=\Lie R$ where $R$ is the affine subgroup of $\Pi$ defined by all
the polynomial equations satisfied by the coefficients of $\dch(\gs)$
plus ``$2\pi i=0$''. Therefore one obtains both
the equality in \eqref{equ:pbdd} and the statement about the completeness
of the standard relations. This concludes the proof of the theorem.
\end{proof}

\begin{rem}

(a) Notice that when $p\ge 11$ the vector space $\ker\beta_p$ contains a subspace
isomorphic to the space
of cusp forms of weight two on $X_1(p)$ which has dimension
$(p-5)(p-7)/24$ (see \cite[Lemma 2.3 \& Theorem 7.8]{G1}).
So it must contain another piece which
has dimension $(p-3)/2$ by \eqref{equ:kerbeta}. What is this missing piece?

(b) When $w>2$ it's not too hard to improve the bound of $d(w,p)$ given
in \cite[5.24]{DG} by the same idea as used in the proof
of the theorem (for example, decrease the bound by $(p^2-1)/24$).
But they are often not the best.
\end{rem}

Turn now to the case of level $N=p^2$. Notice that
$\varphi(p^2)(\varphi(p^2)-2)/4$ is a trivial upper bound of
$\dim \calG_{-2,-2}(\mmu_{p^2})$ (see \cite[Theorem~2.1]{G1}).
If one had any nontrivial upper bound
of $\dim \calG_{-2,-2}(\mmu_{p^2})$ one would have a nontrivial lower
bound of $\dim(\ker \beta_{p^2})$ where $\beta_{p^2}$ is defined
by \eqref{equ:beta}. This would in turn provide an improved
upper bound of $d(2,p^2)$ similar to Theorem~\ref{thm:p}. However,
at present we can only show that
\begin{thm}
If $p\ge 5$ is a prime then $\ker \beta_{p^2}\ne 0$
and $d(2,p^2)<p^2(p-1)^2/4.$
\end{thm}
\begin{proof} Fix a primitive $p^2$th root of unity $\mu$.
Put $e(a)=e_{\mu^a}$ for all integer $a$. Define
$$g_{k,j}=e(pk+j)+e(p^2-pk-j)+e(pj)+e(p^2-pj)$$
for $0\le k<(p-1)/2$, $1\le j\le p-1$, and for $k=(p-1)/2$, $1\le j\le (p-1)/2$.
One only needs to prove the following

\medskip
\noindent
{\bf Claim.} Let $h=(p-3)/2.$ Then one has
\begin{align*}
 \ &\sum_{k=0}^{h} \sum_{l=k}^{h}\sum_{j=2}^{p-2} \{g_{k,1},g_{l,j}\}
+\sum_{k=0}^{h+1} \sum_{j=2}^{h+1} \{g_{k,1},g_{h+1,j}\}
+\sum_{k=0}^{h} \sum_{l=k+1}^{h}\sum_{j=2}^{p-2} \{g_{k,p-1},g_{l,j}\}\\
+&\sum_{k=0}^{h} \sum_{j=2}^{h+1} \{g_{k,p-1},g_{h+1,j}\}
 -\sum_{k=0}^{h} \sum_{l=k}^{h}\sum_{j=2}^{p-2} \{g_{k,j},g_{l,p-1}\}
 -\sum_{k=0}^{h} \sum_{l=k}^{h}\sum_{j=2}^{p-2} \{g_{k,j},g_{l+1,1}\}
 =0.
\end{align*}
There are $h(2h+3)^2=hp^2$ distinct terms on the left, each with coefficient $\pm 1$.

\medskip
The proof of the claim is straight-forward by a little tedious change of indices and regrouping.
$$-\sum_{k=0}^{h} \sum_{l=k}^{h}\sum_{j=2}^{p-2} \{g_{k,j},g_{l+1,1}\}
=\sum_{k=0}^{h} \sum_{l=0}^{k}\sum_{j=2}^{p-2} \{g_{k+1,1},g_{l,j}\}
=\sum_{k=1}^{h+1} \sum_{l=0}^{k-1}\sum_{j=2}^{p-2} \{g_{k,1},g_{l,j}\}.$$
Then the expression in the claim becomes
\begin{align*}
\ & \sum_{k=1}^{h} \sum_{l=0}^{h}\sum_{j=2}^{p-2} \{g_{k,1},g_{l,j}\}
+\sum_{l=0}^{h}\sum_{j=2}^{p-2} \{g_{0,1},g_{l,j}\}
 +\sum_{l=0}^{h}\sum_{j=2}^{p-2} \{g_{h+1,1},g_{l,j}\}
+\sum_{k=0}^{h+1} \sum_{j=2}^{h+1} \{g_{k,1},g_{h+1,j}\} \\
+&\sum_{k=0}^{h} \sum_{l=0}^{h}\sum_{j=2}^{p-2} \{g_{k,p-1},g_{l,j}\}
+\sum_{k=0}^{h} \sum_{j=2}^{h+1} \{g_{k,p-1},g_{h+1,j}\} \\
=\ & \sum_{k=0}^{h+1} \sum_{l=0}^{h}\sum_{j=2}^{p-2} \{g_{k,1},g_{l,j}\}
+\sum_{k=0}^{h+1} \sum_{j=2}^{h+1} \{g_{k,1},g_{h+1,j}\} \\
+&\sum_{k=0}^{h} \sum_{l=0}^{h}\sum_{j=2}^{p-2} \{g_{k,p-1},g_{l,j}\}
+\sum_{k=0}^{h} \sum_{j=2}^{h+1} \{g_{k,p-1},g_{h+1,j}\}
\end{align*}
Let me use $\{a,b\}=\{e(a),e(b)\}$. By definition
\begin{align*}
 \{g_{k,1},g_{l,j}\}=&\{pk+1,pl+j\}
+\{-pk-1,pl+j\}+\{p,pl+j\}+\{-p,pl+j\} \\
+&\{pk+1,-pl-j\}
+\{-pk-1,-pl-j\}+\{p,-pl-j\}+\{-p,-pl-j\} \\
+&\{pk+1,pj\}
+\{-pk-1,pj\}+\{p,pj\}+\{-p,pj\} \\
+&\{pk+1,-pj\}
+\{-pk-1,-pj\}+\{p,-pj\}+\{-p,-pj\}\\
=&\{pk+1,pl+j\}
+\{p(p-k)-1,pl+j\}+\{p,pl+j\}+\{-p,pl+j\} \\
+&\{pk+1,p(p-1-l)+p-j\}
+\{p(p-k)-1,p(p-1-l)+p-j\} \\
\ &\hskip3cm  +\{p,p(p-1-l)+p-j\}+\{-p,p(p-1-l)+p-j\} \\
+&\{pk+1,pj\}
+\{p(p-k)-1,pj\}+\{p,pj\}+\{-p,pj\} \\
+&\{pk+1,p(p-j)\}
+\{p(p-k)-1,p(p-j)\}+\{p,p(p-j)\}+\{-p,p(p-j)\}
\end{align*}
Then
\begin{align*}
\ & \sum_{k=0}^{h+1}  \sum_{l=0}^{h}\sum_{j=2}^{p-2} \{g_{k,1},g_{l,j}\}
=\sum_{k=0}^{h+1} \sum_{l=0}^{h}\sum_{j=2}^{p-2}
\{pk+1,pl+j\}
+\{p(p-k)-1,pl+j\}\\
+&\{pk+1,p(p-1-l)+j\}+\{p(p-k)-1,p(p-1-l)+j\} \\
\ +&\{p,pl+j\}+\{-p,pl+j\} +\{p,p(p-1-l)+j\}+\{-p,p(p-1-l)+j\} \\
+&2\{pk+1,pj\}+2\{p(p-k)-1,pj\}+2\{p,pj\}+2\{-p,pj\} \\
=&\sum_{k=0}^{h+1} \sum_{l=0,l\ne h+1}^{p-1}\sum_{j=2}^{p-2} \{pk+1,pl+j\}
+\sum_{k=h+2}^{p} \sum_{l=0,l\ne h+1}^{p-1}\sum_{j=2}^{p-2} \{pk-1,pl+j\}\\
\ +&\sum_{k=0}^{h+1} \sum_{l=0,l\ne h+1}^{p-1}\sum_{j=2}^{p-2} \Big(\{p,pl+j\}+\{-p,pl+j\}\Big)
+2(h+1)\sum_{k=0}^{h+1}\sum_{j=2}^{p-2} \{pk+1,pj\}\\
+&2(h+1)\sum_{k=h+2}^{p} \sum_{j=2}^{p-2} \{pk-1,pj\}
+2(h+2)(h+1) \sum_{j=2}^{p-2} \Big(\{p,pj\}+\{-p,pj\}\Big)
\end{align*}
Thus
\begin{align*}
\ & \sum_{k=0}^{h+1} \sum_{l=0}^{h}\sum_{j=2}^{p-2} \{g_{k,1},g_{l,j}\}
+\sum_{k=0}^{h+1} \sum_{j=2}^{h+1} \{g_{k,1},g_{h+1,j}\} \\
=&\sum_{k=0}^{h+1} \sum_{l=0}^{p-1}\sum_{j=2}^{p-2} \{pk+1,pl+j\}
+\sum_{k=h+2}^{p} \sum_{l=0}^{p-1}\sum_{j=2}^{p-2} \{pk-1,pl+j\}\\
\ +&\frac{p+1}2\sum_{l=0}^{p-1}\sum_{j=2}^{p-2}\Big(\{p,pl+j\}+\{-p,pl+j\}\Big)
+p\sum_{k=0}^{h+1}\sum_{j=2}^{p-2} \{pk+1,pj\}\\
+&p\sum_{k=h+2}^{p} \sum_{j=2}^{p-2} \{pk-1,pj\}
+\frac{p^2-1}{2} \sum_{j=2}^{p-2} \Big(\{p,pj\}+\{-p,pj\}\Big)
\end{align*}
Similarly,
\begin{align*}
\ & \sum_{k=0}^{h} \sum_{l=0}^{h}\sum_{j=2}^{p-2} \{g_{k,p-1},g_{l,j}\}
+\sum_{k=0}^{h} \sum_{j=2}^{h+1} \{g_{k,p-1},g_{h+1,j}\} \\
=&\sum_{k=1}^{h+1} \sum_{l=0}^{p-1}\sum_{j=2}^{p-2} \{pk-1,pl+j\}
+\sum_{k=h+2}^{p-1} \sum_{l=0}^{p-1}\sum_{j=2}^{p-2} \{pk+1,pl+j\}\\
\ +&\frac{p-1}2\sum_{l=0}^{p-1}\sum_{j=2}^{p-2} \Big(\{p,pl+j\}+\{-p,pl+j\}\Big)
+p\sum_{k=1}^{h+1}\sum_{j=2}^{p-2} \{pk-1,pj\}\\
+&p\sum_{k=h+2}^{p-1} \sum_{j=2}^{p-2} \{pk+1,pj\}
+\frac{(p-1)^2}{2} \sum_{j=2}^{p-2} \Big(\{p,pj\}+\{-p,pj\}\Big)
\end{align*}
Altogether the expression in the claim is reduced to
\begin{align*}
\ &\sum_{k=0}^{p-1} \sum_{l=0}^{p-1}\sum_{j=2}^{p-2} \{pk+1,pl+j\}
 +\sum_{k=1}^{p} \sum_{l=0}^{p-1}\sum_{j=2}^{p-2}\{pk-1,pl+j\} \\
+&p\sum_{l=0}^{p-1}\sum_{j=2}^{p-2} \Big(\{p,pl+j\}+\{-p,pl+j\}\Big)
+p\sum_{k=0}^{p-1}\sum_{j=2}^{p-2} \{pk+1,pj\}\\
+&p\sum_{k=1}^{p}\sum_{j=2}^{p-2} \{pk-1,pj\}
+p(p-1)\sum_{j=2}^{p-2} \Big(\{p,pj\}+\{-p,pj\}\Big)
\end{align*}
This is equal to
\begin{align*}
\sum_{k=0}^{p-1} \sum_{l=0}^{p-1}\sum_{j=2}^{p-2} \Big(&[pk+1,pl+j]-[p(k+l)+j+1,pl+j]+[p(k+l)+j+1,pk+1]\Big)\\
+\sum_{k=1}^p \sum_{l=0}^{p-1}\sum_{j=2}^{p-2} \Big(&[pk-1,pl+j]-[p(k+l)+j-1,pl+j]+[p(k+l)+j-1,pk+1]\Big)\\
+p\sum_{l=0}^{p-1}\sum_{j=2}^{p-2}
\Big(&[p,pl+j]-[p(l+1)+j,pl+j]+[p(l+1)+j,p]\\
+&[-p,pl+j]-[p(l-1)+j,pl+j]+[p(l-1)+j,-p] \Big)\\
+p\sum_{k=0}^{p-1}\sum_{j=2}^{p-2}
\Big(&[pk+1,pj]-[p(j+k)+1,pj]+[p(j+k)+1,pk+1]\Big)\\
+p\sum_{k=1}^{p}\sum_{j=2}^{p-2}
\Big(&[pk-1,pj]-[p(j+k)-1,pj]+[p(j+k)-1,pk-1]\Big)\\
+p(p-1)\sum_{j=2}^{p-2}
\Big(&[p,pj]-[p(j+1),pj]+[p(j+1),p]\\
+&[-p,pj]-[p(j-1),pj]+[p(j-1),-p]\Big)=0
\end{align*}
after many cancelations.

\end{proof}

We have verified the claim in the above proof for $p=5$ and $p=7$ numerically.
In fact, we find that $\dim (\ker \beta_{25})=5$ and $\dim (\ker \beta_{49})=35$.
This implies that $d(2,25)\le 116$ and $d(2,49)\le 449$. Further,
if Grothendieck's transcendence conjecture is true then
the equalities hold. With MAPLE we have verified $d(2,25)\le 116$ and
$d(2,49)\le 449$ by using the standard relations so that these relations
imply all the others.

It is convenient to say a level $N$ is \emph{standard} if either (i)
$N=1, 2$ or $3$, or (ii) $N$ is a prime power $p^n$ ($p\ge 5$).
Otherwise $N$ is called \emph{non-standard}.
We prove in \cite{Zpolyrel} the following
\begin{thm} Suppose $N$ is non-standard. Then the non-standard
relations exist in $\MPV(3,N)$ for $4\le N\le 12$
and in $\MPV(2,N)$ for $10\le N\le 48$.
\end{thm}
Preliminary results suggest that the following
\begin{conj}\label{conj:stand}
If $N$ is a standard level then the standard relations always provide
the sharp bounds of $d(w,N)$, namely, all linear relations can be derived
from the standard ones. If $N$ is a non-standard level then
the bound in \emph{\cite[Cor.\ 5.25]{DG}} is sharp and the non-standard relations
exist in $\MPV(w,N)$ for all $w\ge 3$ (and in $\MPV(2,N)$ if $N\ge 10$).
\end{conj}

\medskip
\noindent
\textbf{Acknowledgements.}
I would like to express my deepest gratitude to Prof. Deligne for
patiently explaining \cite{DG} to me and
suggesting to use octahedral symmetry in level four.
This note was written when I was visiting the Institute for
Advanced Study at Princeton whose financial support and hospitality
is gratefully acknowledged.


\begin{thebibliography}{9}
\bibitem{EZface}
J.~Borwein, P.~Lisonek, and P.~Irvine, \emph{An interface for
evaluation of Euler sums}, available online at
http://oldweb.cecm.sfu.ca/cgi-bin/EZFace/zetaform.cgi


\bibitem{Bourbaki}
N.\ Bourbaki, \emph{Lie Groups and Lie
Algebras Chapters 1-3}, Springer, 2004.

\bibitem{D1}
P.~Deligne, \emph{Le groupe fondamental de la droite projective
  moins trois points}, Galois groups over $\Q$, (Berkeley, CA, 1987),
  Springer, New York, 1989, p.~79--297.

\bibitem{DG}
P.~Deligne and A.~Goncharov, \emph{Groupes fondamentaux motiviques
de Tate mixte}, Ann.\ Scient.\ \'Ec. Norm.
Sup.,  \textbf{38} (1)(2005), 1--56.


\bibitem{G1}
A.~Goncharov, \emph{The dihedral Lie algebras and {G}alois symmetries of
 $\pi_1^{(l)}(\mathbb{P}_1-(\{0,\infty\}\cup\mu_{N}))$}, Duke
  Math.\ J.\ \textbf{110} (3)(2001), 397--487.

\bibitem{GiNac}
J.\ Vollinga, S.\ Weinzierl, \emph{Numerical evaluation of multiple polylogarithms},
Computer Phys. Comm.,  \textbf{167} (3)(2005), 177-194.

\bibitem{Rac}
G.~Racinet, \emph{Doubles m\'elanges des polylogarithmes multiples aux racines de
l'unit\'e}, Publ. Math. Inst. HautesEtudes Sci. \textbf{95}
(2002), 185–231.

\bibitem{Zpolyrel}
J. Zhao, \emph{Standard relations of multiple polylogarithm values of at
roots of unity,} preprint: arxiv:0707.1459v7.

\end{thebibliography}
\end{document}